# Exponentaited Generalized Weibull-Gompertz Distribution

## M.A.EL-Damcese[1], Abdelfattah Mustafa[2] and M.S.Eliwa[2]

[1]Tanta University, Faculty of Science, Mathematics Department, Egypt.

[2]Mansoura University, Faculty of Science, Mathematics Department, Egypt.

**Abstract**

This paper introduces studies on exponentaited generalized Weibull–Gompertz distribution EGWGD($a, b, c, d, \theta$) which generalizes a lot of distributions. Several properties of the EGWGD such as reversed (hazard) function, moments, maximum likelihood estimation, mean residual (past) lifetime, MTTF, MTTR, MTBF, maintainability, availability and order statistics are studied in this paper. A real data set is analyzed and it is observed that the present distribution can provide a better fit than some other very well-known distributions.

***Keywords*:** *Reversed (hazard) function, generalized Weibull-Gompertz distribution, mean residual (past) lifetime, maintainability.*

## 1. Introduction

In analyzing lifetime data one often uses Weibull, Gompertz, Weibull-Gompertz, generalized Gompertz, generalized Weibull-Gompertz distribution and etc. The Weibull distribution has been used in many different fields with many applications see Murthy (2003). The hazard function of the Weibull distribution can only be increasing, decreasing or constant. Thus it cannot be used to model lifetime data with a bathtub shaped hazard function, such as human mortality and machine life cycles. For many years, researchers have been developing various extensions and modified forms of the Weibull distribution, with number of parameters ranging from two to five. Xie et al. (2002) introduced generalization of Weibull which called "the Weibull extension model" and a detailed statistical analysis was given in Tang et al. (2003). The distribution is in fact a generalization of the model studied by Chen (2000).

Also, Lai et al. (2004) proposed a modified Weibull distribution, Bebbington et al. (2007) introduced a flexible Weibull extension and Sarhan et al. (2013) proposed the exponentiated generalized linear exponential distribution which generalizes a large set of distributions including the exponentiated Weibull distribution. These Weibull based parametric models are fit to a breast cancer data set from the National Surgical Adjuvant Breast and Bowel Project.

[*]Corresponding authors:

M.A.EL-Damcese**,** Abdelfattah Mustafa**,**  M.S.Eliwa.

**E-mail**: mseliwa@yahoo.com

On the other hand, The Gompertz distribution is important in describing the pattern of adult deaths see Wetterstrand (1981). For low levels of infant mortality, the Gompertz force of mortality extends to the whole life span see Vaupel (1986) of populations with no observed mortality deceleration. El-Gohary et al. (2013) introduced generalized Gompertz distribution with three parameters. Also, Ali et al. (2014) introduced a new four parameters generalized version of the Gompertz model which is called beta-Gompertz distribution. It includes some well-known lifetime distributions such as beta-exponential and generalized Gompertz distributions as special sub models. Sarhan et al. (2013) proposed exponentiated modified Weibull extension distribution which generalized the distribution mention by El-Gohary et al. (2013). Willekens (2002) provided connections between the Gompertz, the Weibull and other Type I extreme value distributions. Moreover, Nadarajah and Kotz (2005) proposed a generalization of the standard Weibull model with four parameters and this model is called generalized Weibull - Gompertz distribution.

## 2. Exponentiated Generalized Weibull-Gompertz Distribution

The random variable $X$ is said to be has EGWGD if it has the following CDF for $a, b, c, d, \theta > 0$ as follows:

$$F_X(x;a,b,c,d,\theta) = \left[1 - e^{-ax^b(e^{cx^d}-1)}\right]^\theta, \qquad (1)$$

where $b, \theta$ and $d$ are shape parameters, $a$ is scale parameter and $c$ is an acceleration parameter. Exponentiated generalized Weibull-Gompertz distribution with five parameters will denoted by EGWGD$(a, b, c, d, \theta)$. The probability density function $f_X(x;a,b,c,d,\theta)$ of EGWGD$(a, b, c, d, \theta)$ is

$$f_X(x;a,b,c,d,\theta) = ab\theta x^{b-1} e^{-ax^b\left(e^{cx^d}-1\right)+cx^d} \left(1 + \frac{cd}{b}x^d - e^{-cx^d}\right)\left[1 - e^{-ax^b\left(e^{cx^d}-1\right)}\right]^{\theta-1}. \qquad (2)$$

The survival function can be obtained as follows

$$R(x;a,b,c,d,\theta) = 1 - \left[1 - e^{-ax^b\left(e^{cx^d}-1\right)}\right]^\theta \quad ; \ x > 0. \qquad (3)$$

The hazard function $h(x)$ is

$$h(x;a,b,c,d,\theta) = \frac{ab\theta x^{b-1} e^{-ax^b\left(e^{cx^d}-1\right)+cx^d} \left(1 + \frac{cd}{b}x^d - e^{-cx^d}\right)\left[1 - e^{-ax^b(e^{cx^d}-1)}\right]^{\theta-1}}{1 - \left[1 - e^{-ax^b(e^{cx^d}-1)}\right]^\theta}. \qquad (4)$$

Recently, it is observed that the reversed hazard function plays an important role in the reliability analysis; see Gupta and Gupta (2007). The reversed hazard function $r(x)$ of the EGWGD$(a, b, c, d, \theta)$ is

$$r(x;a,b,c,d,\theta) = ab\theta x^{b-1} e^{-ax^b\left(e^{cx^d}-1\right)+cx^d} \left(1 + \frac{cd}{b}x^d - e^{-cx^d}\right)\left[1 - e^{-ax^b(e^{cx^d}-1)}\right]^{-1}. \qquad (5)$$

It is well known that the reversed (hazard) function uniquely determines the corresponding probability density function.

Figures 1 and 3 provide the PDFs of EGWGD$(a, b, c, d, \theta)$ for different parameter values, also figures 2, 4 and 5 provide the failure rate function of EGWGD$(a, b, c, d, \theta)$ for different parameter values, from previous figures it is immediate that the PDFs can be decreasing and unimodal and the hazard function can be increasing, decreasing and bathtub shaped.

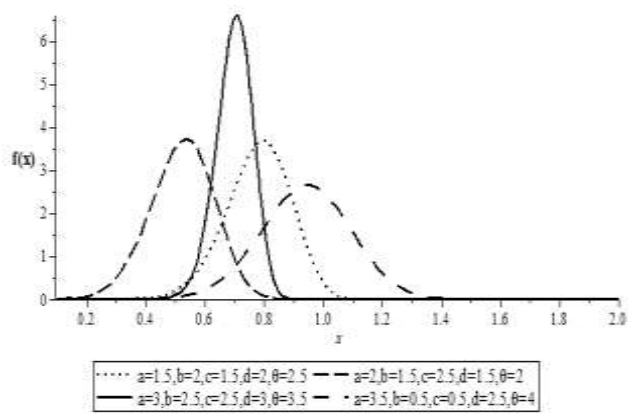

Figure 1. The pdf's of various EGWG distributions

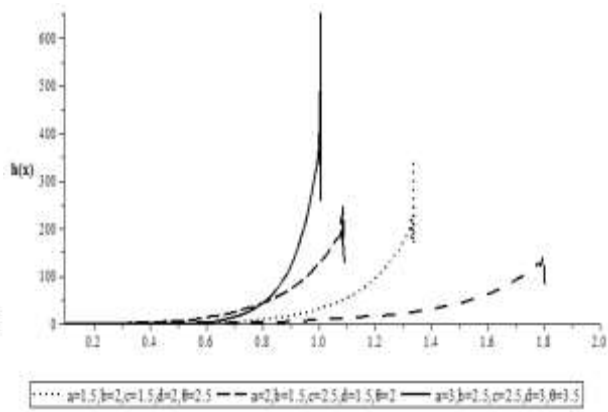

Figure 2. The hazard rate function of various EGWG distributions

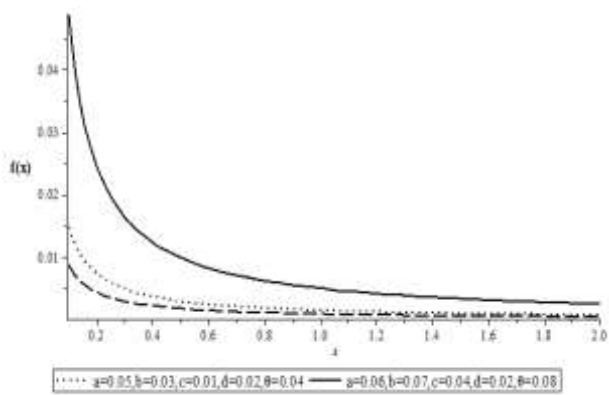

Figure 3. The pdfs of various EGWG distributions

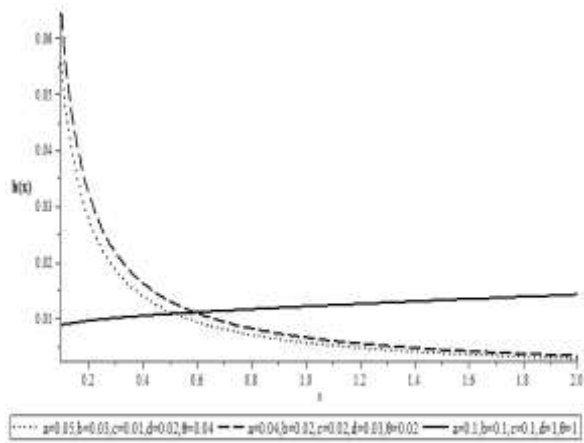

Figure 4. The hazard rate function of various EGWG distributions

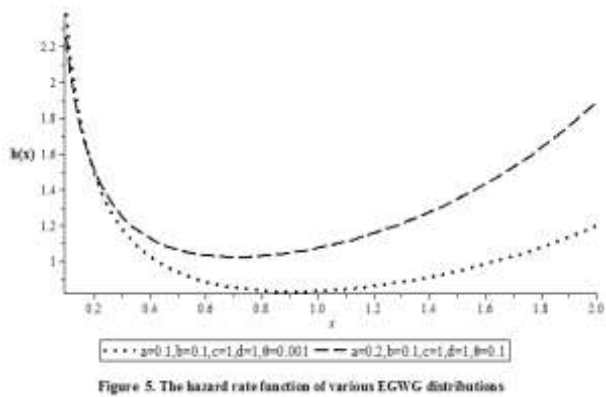

Figure 5. The hazard rate function of various EGWG distributions

**Remark 1.** From EGWGD$(a, b, c, d, \theta)$, the following special cases can be derived:

1. Generalized Weibull - Gompertz distribution GWGD$(a, b, c, d)$, when $\theta = 1$.
2. Gompertz distribution GD $(a, c)$, when $\theta = 1$, $b = 0$ and $d = 1$.
3. Generalized Gompertz distribution GGD $(a, c, \theta)$, when $b = 0$ and $d = 1$.

4. Exponentiated modified Weibull extension, when $b = 0, c = \left(\frac{1}{\alpha}\right)^d, \alpha > 0$.

5. Exponential power distribution $EPD(a, d, c)$, when $\theta = 1$ and $b = 0$.

6. Generalized exponential power distribution $GEPD(a, d, c, \theta)$, when $b = 0$.

7. Weibull extension model of Chen (2000), when $\theta = 1$, $b = 0$ and $c = 1$.

8. Weibull extension model of Xie (2002), when $b = 0$.

9. Exponential distribution $ED(a)$, when $c$ tends to zero, $d = 1, \theta = 1$ and $b = 0$.

10. Generalized exponential distribution $GED(a, \theta)$, when $c$ tends to zero, $d = 1$ and $b = 0$.

## 3. Statistical Properties

### 3.1. The Median and Mode

It is observed as expected that the mean of $EGWGD(a, b, c, d, \theta)$ cannot be obtained in explicit forms. It can be obtained as infinite series expansion so, in general different moments of $EGWGD(a, b, c, d, \theta)$. Also, we cannot get the quartile $x_q$ of $EGWGD(a, b, c, d, \theta)$ in a closed form by using the equation $F_X(x_q; a, b, c, d, \theta) - q = 0$. Thus, by using Equation (1), we find that

$$(x_q)^b e^{c(x_q)^d} = \frac{-1}{a} \ln\left[1 - (q)^{\frac{1}{\theta}}\right], \quad 0 < q < 1. \tag{6}$$

The median $m(X)$ of $EGWGD(a, b, c, d, \theta)$ can be obtained from (6), when $q = 0.5$, as follows

$$(x_{0.5})^b e^{c(x_{0.5})^d} = \frac{-1}{a} \ln\left[1 - (0.5)^{\frac{1}{\theta}}\right]. \tag{7}$$

Moreover, the mode of $EGWGD(a, b, c, d, \theta)$ can be obtained as a solution of the following nonlinear equation.

$$\frac{d}{dx} f_X(x; a, b, c, d, \theta) = 0$$

$$\frac{d}{dx}\left[x^{b-1} e^{-ax^b\left(e^{cx^d}-1\right)+cx^d}\left(1 + \frac{cd}{b}x^d - e^{-cx^d}\right)\left[1 - e^{-ax^b(e^{cx^d}-1)}\right]^{\theta-1}\right] = 0. \tag{8}$$

It is impossible to obtain the explicit solution in general case. It has to be obtained numerically. Explicit forms may be obtained by different special cases.

### 3.2. Moments

The following Lemma 1 gives the $r^{th}$ moment of $EGWGD(a, b, c, d, \theta)$.

**Lemma 1.** If $X$ has $EGWGD(a, b, c, d, \theta)$, the $r^{th}$ moment of $X$, say $\mu'_r$, is given as follows for $a, b, c, d > 0, x > 0$ and $\theta > 1$

Exponentaited Generalized Weibull-Gompertz Distribution

$$\mu'_r = \frac{ab\theta}{d} \sum_{i=0}^{\infty}\sum_{j=0}^{\infty}\sum_{k=0}^{j}\sum_{l=0}^{\infty} \frac{(-1)^{i+j+k} c^l (a(1+i))^j}{j!\, l!\, d\, (ck)^{\frac{r+b(j+1)+ld}{d}}} \binom{j}{k}\binom{\theta-1}{i} \times$$

$$\left( ((1+j)^l - j^l)\Gamma\left(\frac{r+b(j+1)+d(l-1)}{d}+1\right) + \frac{d(1+j)^l}{kb}\Gamma\left(\frac{r+b(j+1)+ld}{d}+1\right)\right). \quad (9)$$

**Proof:**

$$\mu'_r = \int_0^\infty x^r\, f_X(x;a,b,c,d,\theta)\,dx$$

$$= ab\theta \int_0^\infty x^{r+b-1} e^{-ax^b(e^{cx^d}-1)+cx^d}\left(1+\frac{cd}{b}x^d - e^{-cx^d}\right)\left[1-e^{-ax^b(e^{cx^d}-1)}\right]^{\theta-1} dx$$

$$= ab\theta\, I_1 - ab\theta\, I_2 + a\theta cd\, I_3, \quad (10)$$

where

$$I_1 = \int_0^\infty x^{r+b-1} e^{-ax^b(e^{cx^d}-1)+cx^d}\left[1-e^{-ax^b(e^{cx^d}-1)}\right]^{\theta-1} dx,$$

$$I_2 = \int_0^\infty x^{r+b-1} e^{-ax^b(e^{cx^d}-1)}\left[1-e^{-ax^b(e^{cx^d}-1)}\right]^{\theta-1} dx,$$

$$I_3 = \int_0^\infty x^{r+b+d-1} e^{-ax^b(e^{cx^d}-1)+cx^d}\left[1-e^{-ax^b(e^{cx^d}-1)}\right]^{\theta-1} dx,$$

since $0 < \left[1-e^{-ax^b(e^{cx^d}-1)}\right]^{\theta-1} < 1$ for $x > 0$, then by using the binomial series expansion we have

$$\left[1-e^{-ax^b(e^{cx^d}-1)}\right]^{\theta-1} = \sum_{i=0}^{\infty}(-1)^i \binom{\theta-1}{i} e^{-iax^b(e^{cx^d}-1)} \text{ then}$$

$$I_1 = \int_0^\infty x^{r+b-1} e^{cx^d} \sum_{i=0}^{\infty}(-1)^i \binom{\theta-1}{i} e^{-(i+1)ax^b(e^{cx^d}-1)}\,dx,$$

but

$$e^{-(1+i)ax^b(e^{cx^d}-1)} = \sum_{j=0}^{\infty}\sum_{k=0}^{j} \frac{(-1)^{k+j}(a(1+i))^j x^{bj}}{j!}\binom{j}{k} e^{cjx^d} e^{-ckx^d}$$

therefore,

$$I_1 = \sum_{i=0}^{\infty}\sum_{j=0}^{\infty}\sum_{k=0}^{j}\sum_{l=0}^{\infty}\frac{(-1)^{i+j+k}(c(1+j))^l (a(1+i))^j}{j!\, l!}\binom{j}{k}\binom{\theta-1}{i}\int_0^\infty x^{r+b+bj+ld-1} e^{-ckx^d}\,dx$$

Let $\quad y = ckx^d \Rightarrow \quad dx = \frac{1}{d}\left(\frac{y}{ck}\right)^{\frac{1}{d}-1}\left(\frac{1}{ck}\right)dy \quad$ then

$$I_1 = \sum_{i=0}^{\infty}\sum_{j=0}^{\infty}\sum_{k=0}^{j}\sum_{l=0}^{\infty}\frac{(-1)^{i+j+k}(c(1+j))^l(a(1+i))^j}{j!\ l!\ d\ (ck)^{\frac{r+b(j+1)+ld}{d}}}\binom{j}{k}\binom{\theta-1}{i}\int_0^{\infty} y^{\frac{r+b(j+1)+d(l-1)}{d}} e^{-y}\,dy$$

$$=\sum_{i=0}^{\infty}\sum_{j=0}^{\infty}\sum_{k=0}^{j}\sum_{l=0}^{\infty}\frac{(-1)^{i+j+k}(c(1+j))^l(a(1+i))^j}{j!\ l!\ d\ (ck)^{\frac{r+b(j+1)+ld}{d}}}\binom{j}{k}\binom{\theta-1}{i}\Gamma\left(\frac{r+b(j+1)+d(l-1)}{d}+1\right).$$

(11)

Similarly, we find that

$$I_2 = \sum_{i=0}^{\infty}\sum_{j=0}^{\infty}\sum_{k=0}^{j}\sum_{l=0}^{\infty}\frac{(-1)^{i+j+k}(cj)^l(a(1+i))^j}{j!\ l!\ d\ (ck)^{\frac{r+b(j+1)+ld}{d}}}\binom{j}{k}\binom{\theta-1}{i}\Gamma\left(\frac{r+b(j+1)+d(l-1)}{d}+1\right). \quad (12)$$

$$I_3 = \sum_{i=0}^{\infty}\sum_{j=0}^{\infty}\sum_{k=0}^{j}\sum_{l=0}^{\infty}\frac{(-1)^{i+j+k}(c(1+j))^l(a(1+i))^j}{j!\ l!\ d\ (ck)^{\frac{r+b(j+1)+(l+1)d}{d}}}\binom{j}{k}\binom{\theta-1}{i}\Gamma\left(\frac{r+b(j+1)+ld}{d}+1\right). \quad (13)$$

Substituting from (11), (12) and (13) into (10), we get (9). This completes the proof.

## 4. Order Statistics

In fact the order statistics have many applications in reliability and life testing. Let $X_1, X_2, \ldots, X_n$ be a simple random sample from EGWGD with CDF and pdf in Equations (1) and (2) respectively.

Let $X_{(1:n)} \leq X_{(2:n)} \leq \cdots \leq X_{(n:n)}$ denote the order statistics from this sample which $X_{(i:n)}$ denote the lifetime of an $(n-i+1) - out\ of - n$ system which consists of $n$ independent and identically components (iid). Then the pdf of $X_{(i:n)}$, $i = 1,2,3,\ldots,n$ is given by

$$f_{i:n}(x) = \frac{n!}{(i-1)!\ (n-i)!}\ f(x)[F(x)]^{i-1}[1-F(x)]^{n-i}. \quad (14)$$

We defined the first order statistics $X_{(1)} = min(X_1, X_2, \ldots, X_n)$, the last order statistics as $X_{(n)} = max(X_1, X_2, \ldots, X_n)$.

Let $X_1, X_2, \ldots, X_n$ be a simple random sample from EGWGD with CDF and pdf in Equations (1) and (2) respectively. Then the pdf of the first order statistics, $X_{(1)} = X_{1:n}$ is given from equation (14) as follows.

$$f_{1:n}(x) = nab\theta x^{b-1} e^{-ax^b(e^{cx^d}-1)+cx^d}\left(1+\frac{cd}{b}x^d - e^{-cx^d}\right)\left[1-e^{-ax^b(e^{cx^d}-1)}\right]^{\theta-1} \times$$

$$\left(1-\left[1-e^{-ax^b(e^{cx^d}-1)}\right]^{\theta}\right)^{n-1}. \quad (15)$$

Also, the pdf of the last order statistics, $X_{(n)} = X_{n:n}$ is given from Equation (14) as follows.

Exponentaited Generalized Weibull-Gompertz Distribution

$$f_{n:n}(x) = nab\theta x^{b-1} e^{-ax^b(e^{cx^d}-1)+cx^d} \left(1 + \frac{cd}{b}x^d - e^{-cx^d}\right)\left[1 - e^{-ax^b(e^{cx^d}-1)}\right]^{\theta-1} \times$$

$$\left[1 - e^{-ax^b(e^{cx^d}-1)}\right]^{\theta(n-1)}. \qquad (16)$$

## 5. Reliability Analysis

### 5.1 Mean Time to Failure

In order to design and manufacture a maintainable system, it is necessary to predict the mean time to failure (repair), say MTTF (MTTR), for various fault conditions that could occur in the system. This is generally based on past experiences of designers and the expertise available to handle repair work.

The following Lemma 2 obtains the mean time to failure (repair) of the random variable $T$ which has EGWGD$(a,b,c,d,\theta)$.

**Lemma 2.** If $T$ is a random variable has EGWGD$(a,b,c,d,\theta)$, then the mean time to failure (repair) is given as follows for $a, \beta, c, d > 0$ and $\theta > 1$

$$MTTF(MTTR) = \frac{ab\theta}{d} \sum_{i=0}^{\infty}\sum_{j=0}^{\infty}\sum_{k=0}^{j}\sum_{l=0}^{\infty} \frac{(-1)^{i+j+k} c^l (a(1+i))^j}{j!\, l!\, d\, (ck)^{\frac{1+b(j+1)+ld}{d}}} \binom{j}{k}\binom{\theta-1}{i} \times$$

$$\left[((1+j)^l - j^l)\Gamma\left(\frac{1+b(j+1)+d(l-1)}{d}+1\right) + \frac{d(1+j)^l}{kb}\Gamma\left(\frac{1+b(j+1)+ld}{d}+1\right)\right]. \qquad (17)$$

**Proof.** We have

$$\text{MTTF(MTTR)} = \int_0^\infty R(t)\,dt = \int_0^\infty t\, f(t;a,b,c,d,\theta)\,dt = \mu_1',$$

then from Equation (9), when $r = 1$, it is easy to prove Lemma 2.

### 5.2. Availability $A(t)$

It is defined by the probability that the system (component) is successful at time $t$.

**Lemma 3.** If the reliability function for a component is given by $R_T(t) = 1 - F_T(t)$, which $T$ has EGWGD$(a,b,c,d,\theta)$, and the distribution of a repair time density is the PDF $g(t)$ of EGWGD$(a,b,c,d,\theta)$, then the availability $A(t)$ is given as $(t) = 0.5$.

**Proof:**

The proof is simple by using Lemma 2 and $A(t) = MTTF(MTTF + MTTR)^{-1}$.

Moreover, the mean time between failures (MTBF) is an important measure in repairable system (component). This implies that the system (component) has failed and has been repaired. Like MTTF and MTTR, MTBF is the expected value of the random variable time between failures. Mathematically,

$$MTBF = MTTF + MTTR. \tag{18}$$

## 5.3. Maintainability $V(t)$

Let $T$ denote the random variable of the time to repair or the total downtime. If the repair time $T$ has a repair time density function $g(t)$, then the maintainability $V(t)$ is defined as the probability of isolating and repairing a fault in a system within a given time is called maintainability. Or it is the probability that the failed system will be back in service by time $t$. If the repair time $T$ is a random variable has a repair time density function $g(t)$ of EGWGD$(a, b, c, d, \theta)$, then the maintainability $V(t)$ is defined as

$$V(t) = P(T \leq t) = \int_0^t g(s)ds = \left[1 - e^{-at^b\left(e^{ct^d}-1\right)}\right]^\theta \tag{19}$$

## 5.4 The Mean Residual (Past) Life time MRL(MPL) for EGWGD

In reliability theory and survival analysis to study the lifetime characteristics of a live organism there have been defined several measure such as the mean residual lifetime $m(t)$ and the mean past lifetime $p(t)$. Assuming that each component of the system has survived up to time $t$, the survival function of $T_i - t$ given that $T_i > t$; $i = 1, 2, \ldots, n$. This is the corresponding conditional survival function of the components at age $t$.

**Lemma 4.** If $T$ is a random variable has EGWGD$(a, b, c, d, \theta)$, then the mean residual lifetime (MRL), is given as follows for $a, \beta, c, d > 0$ and $\theta > 1$

$$m(t) = \frac{1}{R(t)}[\mu_1' - I(t)] \tag{20}$$

where

$R(t) = 1 - \left[1 - e^{-at^b\left(e^{ct^d}-1\right)}\right]^\theta,$

$\mu_1' = \frac{ab\theta}{d} \sum_{i=0}^{\infty}\sum_{j=0}^{\infty}\sum_{k=0}^{j}\sum_{l=0}^{\infty} \frac{(-1)^{i+j+k} c^l (a(1+i))^j}{j!\, l!\, d\, (ck)^{\frac{1+b(j+1)+ld}{d}}} \binom{j}{k}\binom{\theta-1}{i} \times$

$\left[((1+j)^l - j^l)\Gamma\left(\frac{1+b(j+1)+d(l-1)}{d}+1\right) + \frac{d(1+j)^l}{kb}\Gamma\left(\frac{1+b(j+1)+ld}{d}+1\right)\right],$

$I(t) = t - \sum_{i=0}^{\infty}\sum_{j=0}^{\infty}\sum_{k=0}^{j}\sum_{l=0}^{\infty} \frac{(-1)^{i+j+k+l}(ck)^l(ai)^j}{j!\, l!\, d\, (cj)^{\frac{bj+ld+1}{d}}} \binom{j}{k}\binom{\theta}{i} e^{cjt^d}\left[(cjt^d)^{\frac{bj+ld-d+1}{d}} + \right.$

$\left. \sum_{m=1}^{\frac{bj+ld-d+1}{d}} (-1)^m m! \binom{\frac{bj+ld-d+1}{d}}{m} (cjt^d)^{\frac{bj+ld-d+1}{d}-m}\right].$

**Proof.** By using Equations (1) and (2) we get

Exponentaited Generalized Weibull-Gompertz Distribution

$$m(t) = \frac{1}{R(t)} \int_t^\infty R(x)\, dx = \frac{1}{R(t)} \left( \int_0^\infty x f(x)\, dx - \int_0^t R(x)\, dx \right) = \frac{1}{R(t)} \left( \mu_1' - I(t) \right),$$

where

$$R(t) = 1 - F(t), \quad \mu_1' = \int_0^\infty x f(x)\, dx, \quad I(t) = \int_0^t R(x)\, dx$$

by using Lemma 1, we get

$$\mu_1' = \frac{ab\theta}{d} \sum_{i=0}^\infty \sum_{j=0}^\infty \sum_{k=0}^j \sum_{l=0}^\infty \frac{(-1)^{i+k} c^l (a(1+i))^j}{j!\, l!\, d\, (ck)^{\frac{1+b(j+1)+ld}{d}}} \binom{j}{k} \binom{\theta-1}{i} ((1+j)^l - j^l) \times$$

$$\left[ \Gamma\left( \frac{1 + b(j+1) + d(l-1)}{d} + 1 \right) + \frac{d(1+j)^l}{kb} \Gamma\left( \frac{1 + b(j+1) + ld}{d} + 1 \right) \right].$$

Now, we want to get $I(t)$ as follows.

$$I(t) = \int_0^t R(x)\, dx = t - \int_0^t \left[ 1 - e^{-ax^b(e^{cx^d}-1)} \right]^\theta dx$$

Since

$$\left[ 1 - e^{-ax^b(e^{cx^d}-1)} \right]^\theta = \sum_{i=0}^\infty (-1)^i \binom{\theta}{i} e^{-iax^b(e^{cx^d}-1)}$$

$$= \sum_{i=0}^\infty \sum_{j=0}^\infty \sum_{k=0}^j \sum_{l=0}^\infty \frac{(-1)^{i+j+k+l}(ck)^l(ai)^j}{j!\, l!} \binom{j}{k}\binom{\theta}{i} x^{bj+ld} e^{cjx^d},$$

then

$$I(t) = t - \sum_{i=0}^\infty \sum_{j=0}^\infty \sum_{k=0}^j \sum_{l=0}^\infty \frac{(-1)^{i+j+k+l}(ck)^l(ai)^j}{j!\, l!} \binom{j}{k}\binom{\theta}{i} \int_0^t x^{bj+ld} e^{cjx^d}\, dx.$$

Let $y = cjx^d \quad \Rightarrow \quad dx = \frac{1}{cjd} \left[ \frac{y}{cj} \right]^{\frac{1}{d}-1} dy \quad$ then

$$I(t) = t - \sum_{i=0}^\infty \sum_{j=0}^\infty \sum_{k=0}^j \sum_{l=0}^\infty \frac{(-1)^{i+j+k+l}(ck)^l(ai)^j}{j!\, l!\, d\, (cj)^{\frac{bj+ld+1}{d}}} \binom{j}{k}\binom{\theta}{i} \int_0^{cjt^d} y^{\frac{bj+ld-d+1}{d}} e^y\, dy$$

Since

$$\int_0^u e^y y^n\, dy = e^u \left[ u^n + \sum_{k=1}^n (-1)^k k! \binom{n}{k} u^{n-k} \right]$$

then

$$I(t) = t - \sum_{i=0}^\infty \sum_{j=0}^\infty \sum_{k=0}^j \sum_{l=0}^\infty \frac{(-1)^{i+j+k+l}(ck)^l(ai)^j}{j!\, l!\, d\, (cj)^{\frac{bj+ld+1}{d}}} \binom{j}{k}\binom{\theta}{i} e^{cjt^d} \left[ (cjt^d)^{\frac{bj+ld-d+1}{d}} + \right.$$

$$\left. \sum_{m=1}^{\frac{bj+ld-d+1}{d}} (-1)^m m! \binom{\frac{bj+ld-d+1}{d}}{m} (cjt^d)^{\frac{bj+ld-d+1}{d}-m} \right].$$

this completes the proof.

It is well known that the MRL function $m(t)$ characterizes the distribution function $F(t)$ uniquely, see ( Kotz and Shanbhag (1980)).

The mean past lifetime (MPL) corresponds to the mean time elapsed since the failure of $T_i$ given that $T_i \leq t$. In this case, the random variable of interest is $[t - T_i|T_i \leq t], i = 1,2,3,...,n$. This conditional random variable shows the time elapsed since the failure of $T_i$ given that it failed at or before $t$. The expectation of this random variable gives the mean past lifetime (MPL) $P(t)$.

**Lemma 5.** If $T$ is a random variable has EGWGD$(a, b, c, d, \theta)$, then the mean past lifetime (MPL), is given as follows for $a, b, c, d > 0$ and $\theta > 0$

$$P(t) = \left[1 - e^{-at^b(e^{ct^d}-1)}\right]^{-\theta} \times \left(\sum_{i=0}^{\infty}\sum_{j=0}^{\infty}\sum_{k=0}^{j}\sum_{l=0}^{\infty} \frac{(-1)^{i+j+k+l}(ck)^l(ai)^j}{j!\, l!\, d\, (cj)^{\frac{bj+ld+1}{d}}} \binom{j}{k}\binom{\theta}{i} e^{cjt^d} \times \right.$$

$$\left. \left[(cjt^d)^{\frac{bj+ld-d+1}{d}} + \sum_{m=1}^{\frac{bj+ld-d+1}{d}} (-1)^m m! \binom{\frac{bj+ld-d+1}{d}}{m} (cjt^d)^{\frac{bj+ld-d+1}{d}-m}\right]\right). \quad (21)$$

**Proof:**

$$P(t) = \frac{1}{F(t)} \int_0^t F(x)\, dx = \left[1 - e^{-at^b(e^{ct^d}-1)}\right]^{-\theta} \int_0^t \left[1 - e^{-ax^b(e^{cx^d}-1)}\right]^{\theta} dx. \quad (22)$$

Now, the value of

$$\int_0^t \left[1 - e^{-ax^b(e^{cx^d}-1)}\right]^{\theta} dx = \sum_{i=0}^{\infty}\sum_{j=0}^{\infty}\sum_{k=0}^{j}\sum_{l=0}^{\infty} \frac{(-1)^{i+j+k+l}(ck)^l(ai)^j}{j!\, l!} \binom{j}{k}\binom{\theta}{i} \int_o^t x^{bj+ld} e^{cjx^d} dx,$$

Let $\quad y = cjx^d \quad \Rightarrow \quad dx = \frac{1}{cjd}\left[\frac{y}{cj}\right]^{\frac{1}{d}-1} dy \quad$ then

$$\int_0^t \left[1 - e^{-ax^b(e^{cx^d}-1)}\right]^{\theta} dx = \sum_{i=0}^{\infty}\sum_{j=0}^{\infty}\sum_{k=0}^{j}\sum_{l=0}^{\infty} \frac{(-1)^{i+j+k+l}(ck)^l(ai)^j}{j!\, l!\, d\, (cj)^{\frac{bj+ld+1}{d}}} \binom{j}{k}\binom{\theta}{i} e^{cjt^d} \times$$

$$\left[(cjt^d)^{\frac{bj+ld-d+1}{d}} + \sum_{m=1}^{\frac{bj+ld-d+1}{d}} (-1)^m m! \binom{\frac{bj+ld-d+1}{d}}{m} (cjt^d)^{\frac{bj+ld-d+1}{d}-m}\right], \quad (23)$$

from Equation (23) in Equation (22), we get Equation (21).

The MPL $P(t)$ is also characterizes the underlying distribution uniquely, see Finkelstein (2002).

## 6. Parameters Estimation

### 6.1 Maximum Likelihood Estimates

In this section, we derive the maximum likelihood estimates of the unknown parameters $a, b, c, d$ and $\theta$ of EGWGD$(a, b, c, d, \theta)$, based on a complete sample. Consider a random sample $X_1, X_2, \ldots, X_n$ from EGWGD$(a, b, c, d, \theta)$. The likelihood function of this sample is

$$l = \prod_{i=1}^{n} f(x_i; a, b, c, d, \theta). \tag{24}$$

By substituting from Equation (2) into Equation (24), we get

$$l = \prod_{i=1}^{n} ab\theta x^{b-1} e^{-ax^b(e^{cx^d}-1)+cx^d} \left(1 + \frac{cd}{b} x^d - e^{-cx^d}\right) \left[1 - e^{-ax^b(e^{cx^d}-1)}\right]^{\theta-1}. \tag{25}$$

The log – likelihood function becomes

$$L = n \ln a + n \ln b + n \ln \theta + c \sum_{i=1}^{n}(x_i)^d - a \sum_{i=1}^{n}(x_i)^b \left(e^{c(x_i)^d} - 1\right) + (b-1) \times$$

$$\sum_{i=1}^{n} \ln(x_i) + (\theta - 1) \sum_{i=1}^{n} \ln\left(1 - e^{-a(x_i)^b\left(e^{c(x_i)^d}-1\right)}\right) + \sum_{i=1}^{n} \ln\left(1 + \frac{cd}{b}(x_i)^d - e^{-c(x_i)^d}\right). \tag{26}$$

So, the normal equations are

$$\frac{\partial L}{\partial a} = \frac{n}{\hat{a}} - \sum_{i=1}^{n}(x_i)^{\hat{b}} \left(e^{\hat{c}(x_i)^{\hat{d}}} - 1\right) + (\theta - 1) \sum_{i=1}^{n} \frac{(x_i)^{\hat{b}} \left(e^{\hat{c}(x_i)^{\hat{d}}} - 1\right) e^{-\hat{a}(x_i)^{\hat{b}}\left(e^{\hat{c}(x_i)^{\hat{d}}}-1\right)}}{\left(1 - e^{-\hat{a}(x_i)^{\hat{b}}\left(e^{\hat{c}(x_i)^{\hat{d}}}-1\right)}\right)} = 0. \tag{27}$$

$$\frac{\partial L}{\partial b} = \frac{n}{\hat{b}} - \hat{a} \sum_{i=1}^{n}(x_i)^{\hat{b}} \left(e^{\hat{c}(x_i)^{\hat{d}}} - 1\right) \ln x_i + \sum_{i=1}^{n} \ln x_i - \frac{\hat{c}\hat{d}}{\hat{b}^2} \sum_{i=1}^{n} \frac{(x_i)^{\hat{d}}}{\left(1 + \frac{\hat{c}\hat{d}}{\hat{b}}(x_i)^{\hat{d}} - e^{-c(x_i)^{\hat{d}}}\right)} +$$

$$\hat{a}(\hat{\theta} - 1) \sum_{i=1}^{n}(x_i)^{\hat{b}} \left(e^{\hat{c}(x_i)^{\hat{d}}} - 1\right) e^{-\hat{a}(x_i)^{\hat{b}}\left(e^{\hat{c}(x_i)^{\hat{d}}}-1\right)} \left(1 - e^{-\hat{a}(x_i)^{\hat{b}}\left(e^{\hat{c}(x_i)^{\hat{d}}}-1\right)}\right)^{-1} = 0. \tag{28}$$

$$\frac{\partial L}{\partial c} = \sum_{i=1}^{n}(x_i)^{\hat{d}} - \hat{a} \sum_{i=1}^{n}(x_i)^{\hat{b}+\hat{d}} e^{\hat{c}(x_i)^{\hat{d}}} + \hat{a}(\hat{\theta} - 1) \sum_{i=1}^{n}(x_i)^{\hat{b}+\hat{d}} e^{-\hat{a}(x_i)^{\hat{b}}\left(e^{\hat{c}(x_i)^{\hat{d}}}-1\right)+\hat{c}(x_i)^{\hat{d}}} \times$$

$$\left(1 - e^{-\hat{a}(x_i)^{\hat{b}}\left(e^{\hat{c}(x_i)^{\hat{d}}}-1\right)}\right)^{-1} + \sum_{i=1}^{n} \frac{(x_i)^{\hat{d}} \left(\frac{\hat{d}}{\hat{b}} + e^{-c(x_i)^{\hat{d}}}\right)}{\left(1 + \frac{\hat{c}\hat{d}}{\hat{b}}(x_i)^{\hat{d}} - e^{-c(x_i)^{\hat{d}}}\right)} = 0. \tag{29}$$

$$\frac{\partial L}{\partial d} = \hat{c} \sum_{i=1}^{n}(x_i)^{\hat{d}} \ln x_i - \hat{a}\hat{c} \sum_{i=1}^{n}(x_i)^{\hat{b}+\hat{d}} e^{\hat{c}(x_i)^{\hat{d}}} \ln x_i + \hat{a}\hat{c}(\hat{\theta} - 1) \sum_{i=1}^{n} \frac{(x_i)^{\hat{b}+\hat{d}} e^{\hat{c}(x_i)^{\hat{d}}} \ln x_i}{\left(1 - e^{-\hat{a}(x_i)^{\hat{b}}\left(e^{\hat{c}(x_i)^{\hat{d}}}-1\right)}\right)} +$$

$$\frac{\hat{c}}{\hat{b}} \sum_{i=1}^{n} \frac{(x_i)^{\hat{d}} \left(1 + \hat{d} \ln x_i + \hat{b} e^{-c(x_i)^{\hat{d}}} \ln x_i\right)}{\left(1 + \frac{\hat{c}\hat{d}}{\hat{b}}(x_i)^{\hat{d}} - e^{-c(x_i)^{\hat{d}}}\right)} = 0. \tag{30}$$

The normal equations do not have explicit solution and they have obtained it numerically. The MLE of $\theta$, say $\hat{\theta}(\hat{a}, \hat{b}, \hat{c}, \hat{d})$ can be obtained as

$$\hat{\theta}(\hat{a}, \hat{b}, \hat{c}, \hat{d}) = -n \left( \sum_{i=1}^{n} \ln \left[ 1 - e^{-\hat{a}(x_i)^{\hat{b}} \left( e^{\hat{c}(x_i)^{\hat{d}}} - 1 \right)} \right] \right)^{-1}. \tag{31}$$

So, the MLEs of $\hat{a}, \hat{b}, \hat{c}$ and $\hat{d}$ can be obtained by solving four nonlinear Equations (27) – (30) by using Equation (31).

## 6.2 Asymptotic confidence bounds

In this section, we derive the asymptotic confidence intervals of these parameters when $a, b, c, d > 0$ and $\theta > 0$ as the MLEs of the unknown parameters $a, b, c, d$ can not be obtained in closed forms, by using variance covariance matrix $I_0^{-1}$ see Lawless(2003), where $I_0^{-1}$ is the inverse of the observed information matrix

$$I_0^{-1} = \begin{bmatrix} \frac{-\partial^2 L}{\partial a^2} & \frac{-\partial^2 L}{\partial a \partial b} & \frac{-\partial^2 L}{\partial a \partial c} & \frac{-\partial^2 L}{\partial a \partial d} & \frac{-\partial^2 L}{\partial a \partial \theta} \\ \frac{-\partial^2 L}{\partial b \partial a} & \frac{-\partial^2 L}{\partial b^2} & \frac{-\partial^2 L}{\partial b \partial c} & \frac{-\partial^2 L}{\partial b \partial d} & \frac{-\partial^2 L}{\partial b \partial \theta} \\ \frac{-\partial^2 L}{\partial c \partial a} & \frac{-\partial^2 L}{\partial c \partial b} & \frac{-\partial^2 L}{\partial c^2} & \frac{-\partial^2 L}{\partial c \partial d} & \frac{-\partial^2 L}{\partial c \partial \theta} \\ \frac{-\partial^2 L}{\partial d \partial a} & \frac{-\partial^2 L}{\partial d \partial b} & \frac{-\partial^2 L}{\partial d \partial c} & \frac{-\partial^2 L}{\partial d^2} & \frac{-\partial^2 L}{\partial d \partial \theta} \\ \frac{-\partial^2 L}{\partial \theta \partial a} & \frac{-\partial^2 L}{\partial \theta \partial b} & \frac{-\partial^2 L}{\partial \theta \partial c} & \frac{-\partial^2 L}{\partial \theta \partial d} & \frac{-\partial^2 L}{\partial \theta^2} \end{bmatrix}^{-1}, \tag{32}$$

thus

$$I_0^{-1} = \begin{bmatrix} var(\hat{a}) & cov(\hat{a}, \hat{b}) & cov(\hat{a}, \hat{c}) & cov(\hat{a}, \hat{d}) & cov(\hat{a}, \hat{\theta}) \\ cov(\hat{b}, \hat{a}) & var(\hat{b}) & cov(\hat{b}, \hat{c}) & cov(\hat{b}, \hat{d}) & cov(\hat{b}, \hat{\theta}) \\ cov(\hat{c}, \hat{a}) & cov(\hat{c}, \hat{b}) & var(\hat{c}) & cov(\hat{c}, \hat{d}) & cov(\hat{c}, \hat{\theta}) \\ cov(\hat{d}, \hat{a}) & cov(\hat{d}, \hat{b}) & cov(\hat{d}, c) & var(\hat{d}) & cov(\hat{d}, \hat{\theta}) \\ cov(\hat{\theta}, \hat{a}) & cov(\hat{\theta}, \hat{b}) & cov(\hat{\theta}, \hat{c}) & cov(\hat{\theta}, \hat{d}) & var(\hat{\theta}) \end{bmatrix}. \tag{33}$$

The derivatives in $I_0$ are given as follows:

$$\frac{\partial^2 L}{\partial \theta^2} = \frac{-n}{\theta^2},$$

$$\frac{\partial^2 L}{\partial \theta \partial a} = \sum_{i=1}^{n} \frac{x_i^b \left( e^{cx_i^d} - 1 \right) e^{-ax^b \left( e^{cx_i^d} - 1 \right)}}{1 - e^{-ax_i^b \left( e^{cx_i^d} - 1 \right)}},$$

$$\frac{\partial^2 L}{\partial \theta \partial b} = \sum_{i=1}^{n} \frac{ax_i^b \left( e^{cx_i^d} - 1 \right) e^{-ax_i^b \left( e^{cx_i^d} - 1 \right)} \ln x_i}{\left( 1 - e^{-ax_i^b \left( e^{cx_i^d} - 1 \right)} \right)},$$

$$\frac{\partial^2 L}{\partial \theta \partial c} = \sum_{i=1}^{n} \frac{ax_i^{b+d} e^{-ax_i^b \left( e^{cx_i^d} - 1 \right) + cx_i^d}}{\left( 1 - e^{-ax_i^b \left( e^{cx_i^d} - 1 \right)} \right)},$$

Exponentaited Generalized Weibull-Gompertz Distribution

$$\frac{\partial^2 L}{\partial \theta \partial d} = \sum_{i=1}^{n} \frac{a c x_i^{b+d} e^{-a x_i^b \left(e^{c x_i^d}-1\right)+c x_i^d} \ln x_i}{1 - e^{-a x_i^b \left(e^{c x_i^d}-1\right)}},$$

$$\frac{\partial^2 L}{\partial b^2} = \frac{-n}{b^2} - \sum_{i=1}^{n} x_i^b \left(e^{c x_i^d}-1\right)(\ln x_i)^2 - \frac{cd}{b^2} \sum_{i=1}^{n} \frac{2 b^{-1} x_i^b \left(1 + \frac{cd}{b} x_i^d - e^{-c x_i^d}\right) - \left(\frac{cd}{b^2} x_i^d - e^{-c x_i^d}\right) x_i^d}{\left(1 + \frac{cd}{b} x_i^d - e^{-c x_i^d}\right)^2} +$$

$$a(\theta - 1) \frac{\partial}{\partial b} \sum_{i=1}^{n} \frac{x_i^b \left(e^{c x_i^d}-1\right) e^{-a x_i^b \left(e^{c x_i^d}-1\right)} \ln x_i}{1 - e^{-a x_i^b \left(e^{c x_i^d}-1\right)}},$$

$$\frac{\partial^2 L}{\partial b \partial a} = -\sum_{i=1}^{n} x_i^b \left(e^{c x_i^d}-1\right) \ln x_i + a(\theta - 1) \frac{\partial}{\partial a} \sum_{i=1}^{n} \frac{x_i^b \left(e^{c x_i^d}-1\right) e^{-a x_i^b \left(e^{c x_i^d}-1\right)} \ln x_i}{1 - e^{-a x_i^b \left(e^{c x_i^d}-1\right)}},$$

$$\frac{\partial^2 L}{\partial b \partial c} = -a \sum_{i=1}^{n} x_i^{b+d} e^{c x_i^d} \ln x_i - \frac{d}{b^2} \sum_{i=1}^{n} x_i^d \left(1 + \frac{cd}{b} x_i^d - e^{-c x_i^d}\right) - \frac{c x^{2d} \left(\frac{d}{b} - e^{-c x_i^d}\right)}{\left(1 + \frac{cd}{b} x_i^d - e^{-c x_i^d}\right)^2}$$

$$a(\theta - 1) \frac{\partial}{\partial c} \sum_{i=1}^{n} x_i^b \left(e^{c x_i^d}-1\right) e^{-a x_i^b \left(e^{c x_i^d}-1\right)} \ln x_i \left(1 - e^{-a x_i^b \left(e^{c x_i^d}-1\right)}\right),$$

$$\frac{\partial^2 L}{\partial b \partial d} = -ac \sum_{i=1}^{n} x_i^{b+d} e^{c x_i^d} (\ln x_i)^2 + a(\theta - 1) \frac{\partial}{\partial d} \sum_{i=1}^{n} \frac{x_i^b \left(e^{c x_i^d}-1\right) e^{-a x_i^b \left(e^{c x_i^d}-1\right)} \ln x_i}{\left(1 - e^{-a x_i^b \left(e^{c x_i^d}-1\right)}\right)}$$

$$- \frac{c}{b^2} \sum_{i=1}^{n} \frac{x_i^d \left(1 + \frac{cd}{b} x_i^d - e^{-c x_i^d}\right)(1 + d \ln x_i) - d\, x_i^{2d} \left[\frac{c}{b}(1 + d \ln x_i) + c e^{-c x_i^d} \ln x_i\right]}{\left(1 + \frac{cd}{b} x_i^d - e^{-c x_i^d}\right)^2}$$

$$\frac{\partial^2 L}{\partial a^2} = \frac{-n}{a^2} - (\theta - 1) \sum_{i=1}^{n} \frac{x_i^{2b} \left(e^{c x_i^d}-1\right)^2 e^{-a x_i^b \left(e^{c x_i^d}-1\right)}}{\left(1 - e^{-a x_i^b \left(e^{c x_i^d}-1\right)}\right)^2},$$

$$\frac{\partial^2 L}{\partial a \partial c} = -\sum_{i=1}^{n} x_i^{b+d} e^{c x_i^d} + (\theta - 1) \frac{\partial}{\partial c} \sum_{i=1}^{n} \frac{x_i^b \left(e^{c x_i^d}-1\right) e^{-a x_i^b \left(e^{c x_i^d}-1\right)}}{1 - e^{-a x_i^b \left(e^{c x_i^d}-1\right)}},$$

$$\frac{\partial^2 L}{\partial a \partial d} = -c \sum_{i=1}^{n} x_i^{b+d} e^{c x_i^d} \ln x_i + (\theta - 1) \frac{\partial}{\partial d} \sum_{i=1}^{n} \frac{x_i^b \left(e^{c x_i^d}-1\right) e^{-a x_i^b \left(e^{c x_i^d}-1\right)}}{1 - e^{-a x_i^b \left(e^{c x_i^d}-1\right)}},$$

$$\frac{\partial^2 L}{\partial c^2} = -a \sum_{i=1}^{n} x_i^{b+2d} e^{c x_i^d} - \sum_{i=1}^{n} \frac{x_i^{2d} e^{-c x_i^d} \left(1 + \frac{cd}{b} x_i^d - e^{-c x_i^d}\right) + x_i^{2d} \left(\frac{d}{b} + e^{-c x_i^d}\right)^2}{\left(1 + \frac{cd}{b} x_i^d - e^{-c x_i^d}\right)^2} +$$

$$a(\theta - 1) \frac{\partial}{\partial c} \sum_{i=1}^{n} \frac{x_i^{b+d} e^{-a x_i^b \left(e^{c x_i^d}-1\right)+c x_i^d}}{1 - e^{-a x_i^b \left(e^{c x_i^d}-1\right)}},$$

$$\frac{\partial^2 L}{\partial c \partial d} = \sum_{i=1}^{n} x_i^d \ln x_i - a \sum_{i=1}^{n} x_i^{b+d} e^{cx_i^d} (1 + cx_i^d) \ln x_i$$

$$+ a(\theta - 1) \frac{\partial}{\partial d} \sum_{i=1}^{n} \left( \frac{x_i^{b+d} e^{-ax_i^b(e^{cx_i^d}-1)+cx_i^d}}{1 - e^{-ax_i^b(e^{cx_i^d}-1)}} + \frac{x_i^d \left(\frac{d}{b} + e^{-cx_i^d}\right)}{\left(1 + \frac{cd}{b}x_i^d - e^{-cx_i^d}\right)} \right),$$

$$\frac{\partial^2 L}{\partial d^2} = c \sum_{i=1}^{n} x_i^d (\ln x_i)^2 - ac \sum_{i=1}^{n} x_i^{b+d} e^{cx_i^d} (1 + cx_i^d)(\ln x_i)^2$$

$$+ ac(\theta - 1) \frac{\partial}{\partial d} \sum_{i=1}^{n} \left( \frac{x_i^{b+d} e^{cx_i^d} \ln x_i}{\left(1 - e^{-ax_i^b(e^{cx_i^d}-1)}\right)} + \frac{cx_i^d \left(1 + d \ln x_i + be^{-cx_i^d} \ln x_i\right)}{b\left(1 + \frac{cd}{b}x_i^d - e^{-cx_i^d}\right)} \right),$$

We can derive the $(1 - \delta)100\%$ confidence intervals of the parameters $a, b, c, d, \theta$ by using variance covariance matrix as in the following forms

$\hat{a} \pm Z_{\frac{\delta}{2}} \sqrt{var(\hat{a})}$, $\hat{b} \pm Z_{\frac{\delta}{2}} \sqrt{var(\hat{b})}$, $\hat{c} \pm Z_{\frac{\delta}{2}} \sqrt{var(\hat{c})}$, $\hat{d} \pm Z_{\frac{\delta}{2}} \sqrt{var(\hat{d})}$ and $\hat{\theta} \pm Z_{\frac{\delta}{2}} \sqrt{var(\hat{\theta})}$,

where $Z_{\frac{\delta}{2}}$ is the upper $\left(\frac{\delta}{2}\right)th$ percentile of the standard normal distribution.

## 7. Data Analysis and Discussion

In this section, we present the analysis of a real data set using the EGWGD$(a, b, c, d, \theta)$ model and compare it with the other fitted models like generalized Gompertz distribution GD $(a, c)$, exponential distribution ED$(a)$, Generalized exponential distribution GED$(a, \theta)$, inverse Weibull IW$(\theta)$, generalized inverse Weibull GIW$(\theta, \beta)$, exponentiated inverse Weibull EGIW $(\alpha, \theta, \beta)$. The data have been obtained from Aarset (1987) and it is provided below. It represents the lifetimes of 50 devices.

| 0.1 | 0.2 | 1  | 1  | 1  | 1  | 1  | 2  | 3  | 6  |
|-----|-----|----|----|----|----|----|----|----|----|
| 7   | 11  | 12 | 18 | 18 | 18 | 18 | 18 | 21 | 32 |
| 36  | 40  | 45 | 46 | 47 | 50 | 55 | 60 | 63 | 63 |
| 67  | 67  | 67 | 67 | 72 | 75 | 79 | 82 | 82 | 83 |
| 84  | 84  | 84 | 85 | 85 | 85 | 85 | 85 | 86 | 86 |

The EGWGD model is used to fit this data set. The MLE(s) of the unknown parameter(s), the value of log – likelihood (L), the corresponding Kolmogorove – Smirnove (K-S) ), Akaike information criterion (AIC), correct Akaike information criterion (CAIC) and Bayesian information criterion (BIC) test statistic and its respective p-values for seven different models are given in Table 1

**Table 1.The MLE(s) of the parameter(s),L, AIC, CAIC, BIC, (K-S) values and P-values**

| The Model | MLE(s) | K-S | - L | AIC | CAIC | BIC | P-value |
|---|---|---|---|---|---|---|---|
| ED($a$) | $\hat{a} = 0.022$ | 0.191 | 241.09 | 484.18 | 484.26 | 486.09 | 0.045 |
| GED($a, \theta$) | $\hat{a} = 0.021, \hat{\theta} = 0.902$ | 0.194 | 240.36 | 484.72 | 484.96 | 488.54 | 0.0402 |
| GD ($a, c$) | $\hat{a} = 0.011, \hat{c} = 0.018$ | 0.157 | 235.39 | 474.78 | 475.024 | 478.60 | 0.155 |
| IW($\theta$) | $\hat{\theta} = 0.397$ | 0.435 | 281.07 | 566.14 | 566.39 | 569.96 | $5.9 \times 10^{-9}$ |
| GIW($\theta, \beta$) | $\hat{\theta} = 0.274, \hat{\beta} = 1.273$ | 0.324 | 287.47 | 580.95 | 581.47 | 586.68 | 0.000037 |
| EGIW ($\alpha, \theta, \beta$) | $\hat{\alpha} = 0.75, \hat{\theta} = 0.61$ $\hat{\beta} = 2.142$ | 0.254 | 254.91 | 517.83 | 518.72 | 525.48 | 0.002477 |
| EGWGD (a, b, c, d, θ) | $\hat{a} = 0.000085,$ $\hat{b} = 0.128, \hat{c} = 0.401,$ $\hat{d} = 0.69901, \hat{\theta} = 0.246$ | 0.143 | 224.54 | 467.96 | 469.264 | 477.52 | 0.241 |

Table 1 show that the GD (a,c) and EGWGD(a, b, c, d, θ) are good fit for the given data and in each case we do not reject the hypothesis that the data comes from distribution considered at any usual significance level because the high P − value in each case, but we reject the hypothesis that the data comes from the other distributions because the small P − value in each case. From the P-value it is clear that we reject all the hypotheses when the level of significance is 0.05. But we see that EGWGD(a, b, c, d, θ)model is the best among those distributions because it has the smallest value of (K-S), AIC, CAIC and BIC test. By substituting the MLE of unknown parameters in Equation (13), we get estimation of the variance covariance matrix as

$$I_0^{-1} = \begin{bmatrix} 5.854 \times 10^{-10} & -1.581 \times 10^{-4} & 8.574 \times 10^{-6} & 3.987 \times 10^{-8} & 4.158 \times 10^{-4} \\ -1.581 \times 10^{-4} & 3.101 \times 10^{-5} & -7.004 \times 10^{-4} & -1.012 \times 10^{-5} & 2.175 \times 10^{-4} \\ 8.574 \times 10^{-6} & -7.004 \times 10^{-4} & 3.215 \times 10^{-7} & -5.274 \times 10^{-4} & -2.158 \times 10^{-6} \\ 3.987 \times 10^{-8} & -1.012 \times 10^{-5} & -5.274 \times 10^{-4} & 9.257 \times 10^{-9} & -5.254 \times 10^{-5} \\ 4.158 \times 10^{-4} & 2.175 \times 10^{-4} & -2.158 \times 10^{-6} & -5.254 \times 10^{-5} & 8.154 \times 10^{-5} \end{bmatrix}.$$

The approximate 95% two sided confidence interval of the parameters $a, b, c, d$ and $\theta$ are [0.000037, 0.000132], [0.11708, 0.13891], [0.39988, 0.40211], [0.69881, 0.69918] and [0.2285, 0.26389] respectively.